# Spectral Properties of the Chen–Nagano Gauge on Einstein Manifolds

Sergey Stepanov

ORCID: https://orcid.org/0000-0003-1734-8874

Department of Scientific Information on Fundamental and Applied Mathematics, Russian Institute for Scientific and Technical Information, Russian Academy of Sciences, Moscow, Russia

Department of Mathematics and Data Analysis, Finance University, Moscow, Russia.

E-mail: s.e.stepanov@mail.ru

## Abstract

The stability and deformation theory of Einstein metrics traditionally relies on the classical Berger–Ebin transverse-traceless gauge, which structurally decouples the scalar trace from the divergence-free component of metric perturbations. In the present paper, we introduce a new spectral-geometric framework based on the Chen–Nagano gauge condition. This condition naturally arises from the harmonicity of the identity map and is intrinsically satisfied by the Ricci tensor itself via the contracted second Bianchi identity.

Unlike the classical transverse-traceless framework, the Chen–Nagano gauge preserves a nontrivial interaction between the trace and trace-free sectors of a deformation. We establish a first-order differential relation proving that the divergence of the trace-free part is completely governed by the gradient of the scalar trace. Utilizing commutation formulas on Einstein manifolds, we derive a second-order spectral coupling relation that links the Lichnerowicz Laplacian to a shifted scalar operator.

As a primary geometric consequence, we prove that under suitable spectral pinching assumptions, the Chen–Nagano gauge collapses to the classical transverse-traceless gauge. Specifically, we show that on compact connected negatively curved Einstein manifolds, any volume-preserving Chen–Nagano harmonic deformation whose trace-free component lies below a specific spectral threshold determined by the Einstein constant is necessarily transverse-traceless. Furthermore, we connect this rigidity to the curvature operator of the second kind, establishing explicit lower spectral bounds. Finally, we provide a dynamical interpretation within the Ricci flow framework, demonstrating that the linearized Ricci flow under the Chen–Nagano gauge reduces to a strictly parabolic equation governed by the Lichnerowicz Laplacian, ensuring exponential decay of admissible perturbations.

## 1. Introduction

The stability theory of Einstein metrics occupies a central position in modern differential geometry and geometric analysis (see, for example, [B]; [FP] and [J]). Classical foundations of this theory were established in the works of Berger–Ebin,

Koiso, Fischer–Marsden, and Besse (for details see [B, Chapter 12]), where the second variation of the Hilbert–Einstein functional and the spectral properties of the Lichnerowicz Laplacian were investigated primarily within the framework of transverse-traceless ($TT$) deformations (see [B, Chapter 12]; [K1] and [K2]). In this classical approach, metric deformations are restricted by the conditions

$$\operatorname{tr}_g h = 0, \qquad \delta h = 0$$

for a symmetric 2-tensor $h$ and the divergence operator $\delta,$ which considerably simplify the variational structure of the Einstein equations and the linearization of the Ricci tensor (see [CK, Chapter 3, §2]).

On the other hand, the theory of harmonic symmetric tensors introduced by Chen and Nagano naturally leads to another geometric gauge condition (see [C] and [CN])

$$\delta h = -\frac{1}{2} d(\operatorname{tr}_g h),$$

which will be referred to as the Chen–Nagano gauge. Unlike the classical $TT$-condition, the Chen–Nagano gauge allows nontrivial trace components and therefore defines a geometrically richer class of admissible metric deformations. While classical $TT$-deformations are purely transverse, the Chen–Nagano condition incorporates longitudinal modes naturally coupled to the variation of the trace.

From the geometric point of view, this condition arises from the harmonicity of the identity map between nearby Riemannian metrics and is closely related to harmonic structures on Riemannian manifolds.

Moreover, the Chen–Nagano gauge condition is not merely an auxiliary analytical constraint. It is naturally satisfied by one of the fundamental tensors of Riemannian geometry, namely the Ricci tensor itself, as follows from the contracted second Bianchi identity (see, for example, [CK, p. 73]). This observation indicates that the Chen–Nagano framework is intrinsically related to Einstein geometry and the Ricci flow (see [CK]).

The main purpose of the present paper is to investigate the spectral and geometric structure induced by the Chen–Nagano gauge on Einstein manifolds. We show that this gauge condition produces a nontrivial coupling between the scalar Laplacian acting on the trace part of a deformation and the Lichnerowicz Laplacian acting on its trace-free component. More precisely, if

$$h = h_0 + \frac{1}{n}(\mathrm{tr}_g h)g$$

is the standard pointwise orthogonal decomposition of a symmetric 2-tensor $h$ into its trace-free and pure-trace parts, then the Chen–Nagano condition implies the identity

$$\delta h_0 = -\frac{n-2}{2n} d(\mathrm{tr}_g h).$$

Thus, the divergence of the trace-free component is completely governed by the scalar part of the deformation.

This relation leads to a spectral coupling phenomenon between the scalar Laplacian and the Lichnerowicz Laplacian. In particular, under suitable Einstein and spectral assumptions, the trace part of a Chen–Nagano deformation satisfies a scalar eigenvalue equation whose eigenvalue is explicitly determined by the Einstein constant and the dimension of the manifold. As a consequence, the three-dimensional Einstein case turns out to be exceptional, while in dimensions $n \geq 4$ strong rigidity phenomena arise naturally.

The Chen–Nagano gauge also possesses an important dynamical interpretation. In the sign convention of Besse, the gauge contribution in the linearization of the Ricci tensor cancels identically under the Chen–Nagano condition, and therefore the linearized Ricci flow (see [CK, Chapter 3, §2]) reduces to the strictly parabolic equation

$$\partial_t h = -\Delta_L h.$$

for the Lichnerowicz Laplacian $\Delta_L$. This reduction establishes a direct connection between harmonic symmetric tensors, gauge fixing in geometric analysis, and the dynamical stability theory of Einstein metrics. In particular, scalar-curvature-preserving Chen–Nagano deformations naturally arise in the study of dynamical stability under the Ricci flow (see [CK]). As will be shown below, in the negatively curved Einstein case such deformations necessarily reduce to the classical transverse-traceless gauge.

To put our contribution into perspective, while Fine and Premoselli [FP] address the existence of non-locally symmetric Einstein metrics via topological gluing techniques, and Jäckel [J] investigates their global stability in terms of the large-scale geometry of 3-manifolds, the present work provides an intrinsic analytical description of the space of deformations. We depart from the rigid constraints of the classical transverse-traceless gauge in favor of the Chen–Nagano framework, uncovering a novel spectral coupling phenomenon between the Lichnerowicz Laplacian and the scalar Laplacian. This interaction allows us to establish new spectral and dynamical rigidity bounds that directly govern metric deformations and the linearized Ricci flow on negatively curved Einstein manifolds.

The paper is organized as follows. In Section 2 we recall basic definitions and establish several preliminary identities related to the Chen–Nagano gauge. Section 3 is devoted to the decomposition of symmetric 2-tensors and the induced coupling relations between their trace and trace-free parts. In Section 4 we investigate the corresponding spectral structure and derive rigidity results for negatively curved Einstein manifolds (see [FP] and [J]). Finally, Section 5 discusses geometric consequences and relations to the Ricci flow (see [CK]) and the classical Koiso–Besse stability theory (see [B]).

## 2. Basic Definitions and Preliminaries

Let $(M^n, g)$ be a connected compact Riemannian manifold without boundary. Throughout the paper, all manifolds, metrics, and tensor fields are assumed to be smooth.

Denote by $C^\infty(S^2T^*M)$ the space of smooth symmetric covariant 2-tensor fields on $M^n$. For a tensor $h \in C^\infty(S^2T^*M)$, its local components in a coordinate system $(x^1, \ldots, x^n)$ will be denoted by $h = h_{ij}\, dx^i \otimes dx^j$. The trace of $h$ with respect to the metric $g$ is defined by

$$\mathrm{tr}_g\, h := g^{ij} h_{ij},$$

where $(g^{ij})$ are the contravariant components of the inverse metric tensor $g^{-1}$.

The divergence of $h$ is defined, in the sign convention of Besse, by

$$(\delta h)_j := -\nabla^i h_{ij},$$

where $\nabla^i = g^{ij}\nabla_j$, and $\nabla_j$ denotes the covariant derivative with respect to the coordinate vector field $\frac{\partial}{\partial x^j}$.

The formal adjoint of the divergence operator

$$\delta : C^\infty(S^2 T^* M) \to \Omega^1(M)$$

is the symmetrized covariant derivative

$$\delta^* : \Omega^1(M) \to C^\infty(S^2 T^* M),$$

given by

$$(\delta^* \omega)_{ij} := \frac{1}{2}(\nabla_i \omega_j + \nabla_j \omega_i)$$

where $\omega_j = \omega\left(\frac{\partial}{\partial x^j}\right)$ are the local covariant components of the 1-form $\omega$.

The rough Laplacian acting on symmetric 2-tensors is defined by

$$\nabla^* \nabla h := -\mathrm{tr}_g(\nabla^2 h),$$

while the scalar Laplacian acting on smooth functions is $\Delta := \delta d$. In the convention adopted throughout the paper, the operator $\Delta$ is nonnegative.

We shall also use the Lichnerowicz Laplacian acting on symmetric 2-tensors (see [B, p. 54]):

$$\Delta_L h = \nabla^* \nabla h - 2\mathring{R}(h) + \mathrm{Ric} \circ h + h \circ \mathrm{Ric}\,,$$

where the curvature operator of the second kind is defined by

$$(\mathring{R}(h))_{ij} := R_{ikjl} h^{kl}.$$

The operator $\mathring{R}$ naturally appears as the curvature term in the Lichnerowicz Laplacian. Consequently, the sign of this operator plays a crucial role in rigidity problems for Einstein metrics satisfying

$$\mathrm{Ric} = \lambda g$$

for some constant $\lambda \in \mathbb{R}$. In this case, the Lichnerowicz Laplacian reduces to

$$\Delta_L h = \nabla^* \nabla h - 2\mathring{R}(h) + 2\lambda h.$$

The curvature operator of the second kind plays an important role in rigidity problems for Einstein manifolds and Lichnerowicz-type Laplacians. Classical results show that positivity of this operator leads to strong rigidity phenomena, including spherical space form theorems and vanishing results for harmonic tensors (see [CGT] and [K]). On the other hand, nonpositive and negative curvature operators of the second kind naturally arise on negatively Einstein manifolds and lead to a substantially different spectral picture.

The present work is closely related to earlier investigations devoted to Bochner-type techniques, Lichnerowicz-type Laplacians, and spectral estimates under nonpositive curvature operator assumptions (see [S]). In particular, curvature conditions of the form $\mathring{R} \leq 0$ provide the natural geometric setting for the spectral rigidity and transverse-traceless reduction phenomena established in the present paper.

A central role in the present work is played by the Chen–Nagano gauge condition

$$\delta h = -\frac{1}{2} d(\mathrm{tr}_g h),$$

which naturally arises from the harmonicity of the identity map (see [C]; [CN])

$$\mathrm{id}_M : (M^n, g) \to (M^n, g + h)$$

and may be interpreted as a geometric gauge condition associated with the contracted Bianchi identity. Symmetric 2-tensors satisfying the Chen–Nagano condition will be called Chen–Nagano harmonic tensors.

Unlike the classical transverse-traceless conditions

$$\operatorname{tr}_g h = 0, \qquad \delta h = 0,$$

the Chen–Nagano gauge allows nontrivial trace components and therefore defines a broader class of admissible metric deformations.

It is well known that the space of transverse-traceless tensors on a compact manifold is infinite-dimensional (see [B, §4.57 and §4.61]). Since every transverse-traceless tensor satisfies the Chen–Nagano condition, the space of Chen–Nagano harmonic symmetric tensors contain the classical $TT$-subspace and therefore is also infinite-dimensional.

Every symmetric 2-tensor admits the pointwise orthogonal decomposition

$$h = h_0 + \frac{1}{n}(\operatorname{tr}_g h)g,$$

where

$$\operatorname{tr}_g h_0 = 0.$$

The Chen–Nagano condition imposes a nontrivial relation between the divergence of the trace-free component $h_0$ and the scalar part $\operatorname{tr}_g h$. This relation will serve as the starting point for the spectral analysis developed in the subsequent sections.

## 3. The Chen–Nagano Gauge and Spectral Coupling

In this section, we investigate the interaction between the Chen–Nagano gauge condition and the standard orthogonal decomposition of symmetric 2-tensors into their trace-free and pure-trace parts. We shall show that the Chen–Nagano condition induces a nontrivial coupling between these components and therefore leads naturally to a spectral interaction between the scalar Laplacian and the Lichnerowicz Laplacian.

Let $h \in C^\infty(S^2T^*M)$ be a symmetric 2-tensor. Recall that every such tensor admits the pointwise orthogonal decomposition

$$h = h_0 + \frac{1}{n}(\operatorname{tr}_g h)g,$$

where

$$\operatorname{tr}_g h_0 = 0.$$

The following result describes the first-order structure induced by the Chen–Nagano gauge.

**Proposition 3.1.** *Let* $(M^n, g)$ *be a Riemannian manifold and let*

$$h = h_0 + \frac{1}{n}(tr_g h)g$$

*be the pointwise orthogonal decomposition of a symmetric* $2$*-tensor into its trace-free and pure-trace parts. If h satisfies the Chen–Nagano gauge condition*

$$\delta h = -\frac{1}{2} d(tr_g h),$$

*then*

$$\delta h_0 = -\frac{n-2}{2n} d(tr_g h).$$

**Proof.** Using the decomposition of $h$, we obtain

$$\delta h = \delta h_0 + \frac{1}{n}\delta((\operatorname{tr}_g h)g).$$

Since the Levi-Civita connection is metric-compatible, $\nabla g = 0$, it follows from the definition of the divergence operator that

$$\delta(fg) = -df$$

for every smooth function $f \in C^\infty(M)$. Therefore,

$$\delta((\operatorname{tr}_g h)g) = -d(\operatorname{tr}_g h).$$

Substituting this identity into the previous formula yields

$$\delta h = \delta h_0 - \frac{1}{n} d(\operatorname{tr}_g h).$$

Since $h$satisfies the Chen–Nagano condition,

$$\delta h = -\frac{1}{2} d(\operatorname{tr}_g h),$$

we obtain

$$\delta h_0 - \frac{1}{n} d(\mathrm{tr}_g h) = -\frac{1}{2} d(\mathrm{tr}_g h).$$

Hence,

$$\delta h_0 = \left(\frac{1}{n} - \frac{1}{2}\right) d(\mathrm{tr}_g h),$$

which gives

$$\delta h_0 = -\frac{n-2}{2n} d(\mathrm{tr}_g h).$$

The proof is complete. □

Proposition 3.1 shows that, unlike the classical transverse-traceless theory, the divergence of the trace-free component $h_0$ does not vanish in general. Instead, it is completely determined by the scalar part of the deformation. Thus, within the Chen–Nagano framework, the trace-free and scalar components are intrinsically coupled.

In particular, if $n = 2$, then Proposition 3.1 implies

$$\delta h_0 = 0.$$

Hence, in dimension two, the trace-free component of every Chen–Nagano harmonic tensor is automatically divergence-free.

For dimensions $n \geq 3$, the divergence of $h_0$ is entirely controlled by the gradient of the trace. This coupling phenomenon serves as the starting point for the spectral analysis developed below.

To understand how the first-order coupling established in Proposition 3.1 gives rise to a spectral interaction between the Lichnerowicz Laplacian and the scalar Laplacian, we now investigate the divergence of $\Delta_L h_0$.

Throughout this section, we assume that $(M^n, g)$ is an Einstein manifold satisfying

$\mathrm{Ric} = \lambda g$ for some constant $\lambda \in \mathbb{R}$. Recall the classical commutation relation between the divergence operator and the Lichnerowicz Laplacian acting on symmetric 2-tensors (see [B, Chapter 4]):

$$\delta(\Delta_L h) = \Delta_H(\delta h),$$

where

$$\Delta_H = d\delta + \delta d$$

is the Hodge–de Rham Laplacian acting on 1-forms.

Applying this identity to the trace-free component $h_0$ and using Proposition 3.1, we obtain the following differential relation.

**Theorem 3.2.** *Let* $(M^n, g)$ *be a compact Einstein manifold satisfying* $Ric = \lambda g$. *Let*

$$h = h_0 + \frac{1}{n}(tr_{\,g}\, h)g$$

*be a Chen–Nagano harmonic tensor. Then*

$$\delta(\Delta_L h_0) = -\frac{n-2}{2n} d\big(\Delta(tr_{\,g}\, h) - 2\lambda\, tr_{\,g}\, h\big).$$

**Proof.** Applying the commutation relation

$$\delta(\Delta_L h) = \Delta_H(\delta h) - 2\lambda\, \delta h$$

to the trace-free component $h_0$, we obtain

$$\delta(\Delta_L h_0) = \Delta_H(\delta h_0) - 2\lambda\, \delta h_0.$$

By Proposition 3.1,

$$\delta h_0 = -\frac{n-2}{2n} d(\mathrm{tr}_{\,g}\, h).$$

Therefore,

$$\delta(\Delta_L h_0) = \Delta_H\left(-\frac{n-2}{2n} d(\mathrm{tr}_{\,g}\, h)\right) - 2\lambda\left(-\frac{n-2}{2n} d(\mathrm{tr}_{\,g}\, h)\right).$$

Since the Hodge–de Rham Laplacian commutes with the exterior derivative,

$$\Delta_H d = d\Delta,$$

we obtain

$$\delta(\Delta_L h_0) = -\frac{n-2}{2n} d\left(\Delta\big(\mathrm{tr}_g h\big)\right) + \frac{n-2}{2n} d(2\lambda\, \mathrm{tr}_{\,g}\, h).$$

Since $\lambda$ is constant, this gives

$$\delta(\Delta_L h_0) = -\frac{n-2}{2n} d\big(\Delta(\operatorname{tr}_g h) - 2\lambda \operatorname{tr}_g h\big).$$

The proof is complete. □

Theorem 3.2 reveals the precise mechanism underlying the spectral interaction induced by the Chen–Nagano gauge. It shows that the divergence of $\Delta_L h_0$ is completely governed by the shifted scalar Laplacian acting on the trace part of the deformation. Thus, within the Chen–Nagano framework, the scalar and trace-free components of a metric deformation are no longer spectrally independent. Instead, they are linked through a differential coupling relation connecting the scalar Laplacian and the Lichnerowicz Laplacian. This interaction will serve as the analytic foundation for the rigidity results established in the next section.

## 4. Spectral Rigidity and the Curvature Operator of the Second Kind

In this section, we derive rigidity consequences of the spectral coupling established in Section 3. We show that, under suitable spectral assumptions on the trace-free component $h_0$, the Chen–Nagano gauge collapses to the classical transverse-traceless gauge. We further relate this phenomenon to the curvature operator of the second kind.

Throughout this section, let $(M^n, g)$ be a compact connected Einstein manifold satisfying $\mathrm{Ric} = \lambda g,\ \lambda < 0$. Assume that

$$h = h_0 + \frac{1}{n}(\operatorname{tr}_g h) g$$

is a Chen–Nagano harmonic tensor.

We begin with the principal rigidity theorem.

**Theorem 4.1.** *Let $(M^n, g)$, $n \geq 3$, be a compact connected Einstein manifold satisfying*

$$Ric = \lambda g, \qquad \lambda < 0.$$

*Assume that $h$ is a Chen–Nagano harmonic tensor and that its trace-free component $h_0$ satisfies the eigenvalue equation*

$$\Delta_L h_0 = \mu h_0$$

*for some* $\mu \in \mathbb{R}$*. If*

$$\mu < -2\lambda,$$

*then the scalar trace of h is constant:*

$$\operatorname{tr}_g h = C_0$$

*for some constant* $C_0 \in \mathbb{R}$*. Consequently,* $\delta h_0 = 0$.

*Furthermore, if h satisfies the first-order volume-preserving condition*

$$\int_M \operatorname{tr}_g h \, dV_g = 0,$$

*then* $\operatorname{tr}_g h \equiv 0$. *Hence*

$$h = h_0, \qquad \delta h = 0, \qquad \operatorname{tr}_g h = 0,$$

*so h is a classical transverse-traceless tensor.*

**Proof.** Taking the divergence of the eigenvalue equation

$$\Delta_L h_0 = \mu h_0,$$

we obtain

$$\delta(\Delta_L h_0) = \mu \, \delta h_0.$$

Using Proposition 3.1,

$$\delta h_0 = \frac{2-n}{2n} d(\operatorname{tr}_g h),$$

we get

$$\delta(\Delta_L h_0) = \mu \frac{2-n}{2n} d(\operatorname{tr}_g h).$$

On the other hand, by Theorem 3.2,

$$\delta(\Delta_L h_0) = -\frac{n-2}{2n} d(\Delta(\operatorname{tr}_g h) - 2\lambda \operatorname{tr}_g h).$$

Since $n \geq 3$, we may cancel the nonzero factor $\frac{2-n}{2n}$. Therefore,

$$d(\Delta(\operatorname{tr}_g h) - 2\lambda \operatorname{tr}_g h) = \mu \, d(\operatorname{tr}_g h),$$

or equivalently,

$$d(\Delta(\operatorname{tr}_g h) - (2\lambda + \mu)\operatorname{tr}_g h) = 0.$$

Since $M^n$ is connected and compact, it follows that

$$\Delta(\operatorname{tr}_g h) - (2\lambda + \mu)\operatorname{tr}_g h = C$$

for some constant $C \in \mathbb{R}$.

Set $u = \operatorname{tr}_g h$. Since

$$\mu < -2\lambda,$$

we have

$$2\lambda + \mu < 0.$$

Let

$$a = -(2\lambda + \mu) > 0.$$

Then the above equation becomes

$$\Delta u + au = C.$$

Define

$$\bar{u} = \frac{1}{\mathrm{Vol}\,(M)} \int_M u \, d\, V_g$$

and let $v = u - \bar{u}$. Then

$$\int_M v \, d\, V_g = 0$$

and $\Delta v + av = 0$. Multiplying by $v$and integrating over $M^n$, we obtain

$$\int_M v\Delta v \, d\, V_g + a \int_M v^2 \, dV_g = 0.$$

Integrating by parts gives

$$\int_M |\, dv \,|^2 \, dV_g + a \int_M v^2 \, dV_g = 0.$$

Since $a > 0$, both terms are nonnegative. Hence,

$$v \equiv 0.$$

Therefore, $u = \operatorname{tr}_g h$ is constant:

$$\operatorname{tr}_g h = C_0.$$

Using Proposition 3.1 again, we conclude that $\delta h_0 = 0$.

Finally, if

$$\int_M \operatorname{tr}_g h \, dV_g = 0,$$

then the constant $C_0$vanishes identically. Hence $\operatorname{tr}_g h = 0$. Consequently,

$$h = h_0, \qquad \delta h = 0, \qquad \operatorname{tr}_g h = 0.$$

Thus $h$ is a transverse-traceless tensor. The proof is complete. □

To understand the geometric meaning of the spectral condition appearing in Theorem 4.1, we now relate the eigenvalue $\mu$ to the curvature operator of the second kind.

Recall that on an Einstein manifold satisfying $\mathrm{Ric} = \lambda g$, the Lichnerowicz Laplacian acting on trace-free symmetric 2-tensors admits the Weitzenböck representation

$$\Delta_L h_0 = \nabla^* \nabla h_0 - 2\mathring{R}(h_0) + 2\lambda h_0.$$

Assume that

$$\Delta_L h_0 = \mu h_0.$$

Taking the $L^2$-inner product with $h_0$ and integrating over $M^n$, we obtain

$$\mu \int_M | h_0 |^2 \, dV_g = \int_M | \nabla h_0 |^2 \, dV_g - 2 \int_M \langle \mathring{R}(h_0), h_0 \rangle \, dV_g + 2\lambda \int_M | h_0 |^2 \, dV_g.$$

We now use the classical Bochner–Koiso identity for trace-free symmetric 2-tensors on Einstein manifolds:

$$\int_M | \nabla h_0 |^2 \, dV_g = \int_M | \delta h_0 |^2 \, dV_g + \int_M \langle \mathring{R}(h_0), h_0 \rangle \, dV_g + \lambda \int_M | h_0 |^2 \, dV_g.$$

Substituting this identity into the previous formula yields

$$\mu \int_M | h_0 |^2 \, dV_g = \int_M | \delta h_0 |^2 \, dV_g - \int_M \langle \mathring{R}(h_0), h_0 \rangle \, dV_g + 3\lambda \int_M | h_0 |^2 \, dV_g.$$

Equivalently,

$$(\mu - 3\lambda) \int_M | h_0 |^2 \, dV_g = \int_M | \delta h_0 |^2 \, dV_g - \int_M \langle \mathring{R}(h_0), h_0 \rangle \, dV_g.$$

This identity immediately leads to the following estimate.

**Corollary 4.2.** *Let $(M^n, g)$ be a compact connected Einstein manifold satisfying*

$$Ric = \lambda g, \qquad \lambda < 0.$$

*Assume that the curvature operator of the second kind satisfies*

$$\mathring{R} \leq 0$$

*on trace-free symmetric 2-tensors. If*

$$\Delta_L h_0 = \mu h_0$$

*for a nontrivial trace-free tensor $h_0$, then*

$$\mu \geq 3\lambda.$$

**Proof.** Since $\mathring{R} \leq 0$, we have

$$- \int_M \langle \mathring{R}(h_0), h_0 \rangle \, dV_g \geq 0.$$

Moreover,

$$\int_M | \delta h_0 |^2 \, dV_g \geq 0.$$

Hence the right-hand side of

$$(\mu - 3\lambda) \int_M | h_0 |^2 \, dV_g = \int_M | \delta h_0 |^2 \, dV_g - \int_M \langle \mathring{R}(h_0), h_0 \rangle \, dV_g$$

is nonnegative.

Since $h_0 \not\equiv 0$,

$$\int_M | h_0 |^2 \, dV_g > 0.$$

Therefore, $\mu - 3\lambda \geq 0$, which proves

$$\mu \geq 3\lambda.$$

The proof is complete. □

It is worth noting that, for negatively Einstein manifolds ($\lambda < 0$), Corollary 4.2 shows that, under the curvature condition $\mathring{R} \leq 0$, the spectrum of the Lichnerowicz Laplacian on trace-free Chen–Nagano harmonic eigentensors is bounded below by $3\lambda$.

On the other hand, Theorem 4.1 implies that every eigendeformation satisfying

$$\mu < -2\lambda$$

necessarily has constant trace and therefore reduces, under the first-order volume-preserving condition, to the classical transverse-traceless gauge.

Consequently, the possible spectral range for nontrivial Chen–Nagano harmonic deformations with nonconstant trace is confined to the interval $[3\lambda, -2\lambda)$. This demonstrates that the Chen–Nagano gauge possesses a rigid spectral structure on negatively Einstein manifolds: once the eigenvalue crosses the threshold $-2\lambda$, the deformation is forced into the classical transverse-traceless regime.

## 5. Dynamical Stability and the Linearized Ricci Flow

In this final section, we discuss the dynamical implications of the spectral rigidity phenomena established in the previous sections. We show that, under the Chen–Nagano gauge condition, the linearized Ricci flow (see [CK, Chapter 3, §2]) reduces to a purely parabolic evolution equation governed by the Lichnerowicz Laplacian. Combined with the spectral coupling results of Sections 3 and 4, this leads naturally to dynamical rigidity and exponential decay of admissible perturbations.

Recall that the Ricci flow (see [CK]) is defined by

$$\partial_t g = -2\mathrm{Ric}\,(g).$$

A remarkable feature of the Chen–Nagano gauge is its compatibility with the linearization of the Ricci tensor. In the sign convention of Besse, the linearization of the Ricci tensor is given by (see [B, p. 63])

$$(\mathrm{Ric}\,)'(h) = \frac{1}{2}(\Delta_L h - 2\delta^*\delta h - dd(\mathrm{tr}_g\, h)).$$

Therefore, under the Chen–Nagano gauge condition

$$\delta h = -\frac{1}{2} d(\mathrm{tr}_g\, h),$$

the gauge contribution vanishes identically:

$$-2\delta^*\delta h - dd(\mathrm{tr}_g\, h) = 0.$$

Consequently, the linearized Ricci tensor reduces to

$$(\mathrm{Ric}\,)'(h) = \frac{1}{2}\Delta_L h.$$

Thus, the linearization of the Ricci flow (see [CK, Chapter 3]) at an Einstein metric takes the form

$$\partial_t h = -\Delta_L h.$$

Hence, under the Chen–Nagano gauge, the linearized Ricci flow (see [CK, Chapter 3, §2]) becomes a strictly parabolic equation whose dynamics are completely determined by the spectrum of the Lichnerowicz Laplacian.

We now combine this observation with the rigidity theory developed in Section 4.

**Theorem 5.1.** *Let $(M^n, g)$, $n \geq 3$, be a compact connected Einstein manifold satisfying*

$$Ric = \lambda g, \qquad \lambda < 0.$$

*Assume that his a volume-preserving Chen–Nagano harmonic perturbation whose trace-free component satisfies the eigenvalue equation*

$$\Delta_L h_0 = \mu h_0$$

*for some eigenvalue*

$$0 < \mu < -2\lambda.$$

*Then $h$ reduces to a classical transverse-traceless perturbation and decays exponentially under the linearized Ricci flow.*

**Proof.** Since

$$\mu < -2\lambda,$$

Theorem 4.1 implies that

$$\mathrm{tr}_g\, h \equiv 0.$$

Therefore,

$$h = h_0, \qquad \delta h = 0,$$

and hence $h$is a classical transverse-traceless tensor.

Under the linearized Ricci flow (see [CK, Chapter 3, §2])

$$\partial_t h = -\Delta_L h,$$

the eigentensor condition

$$\Delta_L h = \mu h$$

implies $\partial_t h = -\mu h$. Therefore, the solution is given explicitly by

$$h(t) = e^{-\mu t} h(0).$$

Since $\mu > 0$, it follows that

$$\lim_{t \to \infty} h(t) = 0$$

exponentially fast. More precisely,

$$\| h(t) \|_{L^2} = e^{-\mu t} \| h(0) \|_{L^2}.$$

The proof is complete. □

Theorem 5.1 shows that the Chen–Nagano gauge exhibits a strong dynamical rigidity phenomenon on negatively Einstein manifolds. Under the spectral pinching condition

$$0 < \mu < -2\lambda,$$

every admissible Chen–Nagano harmonic perturbation is forced into the classical transverse-traceless regime and undergoes exponential decay under the linearized Ricci flow (see [CK, Chapter 3, §2]).

The results obtained in the present paper demonstrate that the Chen–Nagano gauge provides a natural geometric framework in which the scalar and trace-free components of Einstein metric deformations remain spectrally coupled. This coupling mechanism leads to rigidity phenomena, spectral restrictions for the Lichnerowicz Laplacian, and dynamical stability results for the linearized Ricci flow (see [CK, Chapter 3, §2]).

In this sense, the Chen–Nagano approach extends the classical Berger–Ebin transverse-traceless theory by replacing the complete decoupling of the scalar and divergence-free components with a nontrivial spectral interaction governed by the geometry of the Einstein manifold and the curvature operator of the second kind.

The appearance of the Ricci tensor in the Chen–Nagano framework is not accidental. Indeed, the contracted second Bianchi identity implies

$$\delta \mathrm{Ric} = -\frac{1}{2} ds.$$

Since $\mathrm{tr}_g(\mathrm{Ric}) = s$, the Ricci tensor satisfies precisely the Chen–Nagano gauge condition

$$\delta \mathrm{Ric} = -\frac{1}{2} d(\mathrm{tr}_g \mathrm{Ric})$$

on every Riemannian manifold (see [CK, p. 73]).

Thus, the Chen–Nagano condition is not merely an auxiliary analytical constraint, but is intrinsically encoded in the geometry of the Ricci tensor itself. This observation explains the natural compatibility of the Chen–Nagano framework with Einstein geometry and the Ricci flow.

Furthermore, if the scalar curvature of a Riemannian manifold $(M^n, g)$ vanishes, $s = 0$, then $\operatorname{tr}_g(\mathrm{Ric}) = 0$ and $\delta\mathrm{Ric} = 0$. Consequently, in the scalar-flat case, the Ricci tensor becomes a genuine transverse-traceless tensor.

**Conclusion**

In the present paper, we introduced a new spectral-geometric approach to the study of Einstein metric deformations based on the Chen–Nagano gauge condition

$$\delta h = -\frac{1}{2} d(\operatorname{tr}_g h)$$

for a symmetric 2-tensor $h$, which naturally arises from the harmonicity of the identity map

$$\mathrm{id}_M (M^n, g) \to (M^n, g + h).$$

Unlike the classical Berger–Ebin transverse-traceless decomposition, where the trace-free and divergence-free conditions are imposed independently,

$$\operatorname{tr}_g h = 0, \qquad \delta h = 0,$$

the Chen–Nagano framework preserves a nontrivial interaction between the scalar and trace-free components of a deformation. More precisely, the Chen–Nagano condition produces the differential coupling relation

$$\delta h_0 = \frac{2-n}{2n} d(\operatorname{tr}_g h),$$

which shows that the divergence of the trace-free component $h_0$ of $h$ is completely governed by the scalar trace.

From the geometric point of view, this makes the Chen–Nagano gauge substantially richer and technically more delicate than the classical transverse-traceless theory. In the Berger–Ebin framework, the scalar and trace-free sectors decouple completely. By contrast, the Chen–Nagano condition forces these components to interact through a coupled spectral structure involving simultaneously the scalar Laplacian and the Lichnerowicz Laplacian.

Using the commutation formulas for the divergence operator and the Lichnerowicz Laplacian on Einstein manifolds, we established a second-order spectral coupling relation connecting the operator $\Delta_L$ with the shifted scalar operator $\Delta - (2\lambda + \mu)$. This interaction led to a series of rigidity results showing that, under suitable spectral pinching assumptions, the Chen–Nagano gauge collapses to the classical transverse-traceless gauge. In particular, we proved that negatively Einstein manifolds satisfying $\mu < -2\lambda$ admit no nontrivial volume-preserving Chen–Nagano harmonic deformations outside the classical $TT$-regime.

A further geometric aspect of the theory arises from the curvature operator of the second kind. By combining Bochner–Koiso identities with the Weitzenböck formula for the Lichnerowicz Laplacian, we derived explicit lower spectral bounds controlled by the sign of the curvature operator. This demonstrates that the curvature operator of the second kind governs not only rigidity properties of Einstein metrics, but also the spectral behavior of Chen–Nagano harmonic deformations.

Finally, the compatibility of the Chen–Nagano gauge with the linearization of the Ricci tensor allowed us to connect the entire theory with the dynamical stability problem for the Ricci flow (see [CK]). Under the Chen–Nagano condition, the linearized Ricci flow (see [CK, Chapter 3, §2]) reduces to a purely parabolic equation governed by the Lichnerowicz Laplacian. As a consequence, admissible eigendeformations satisfying $0 < \mu < -2\lambda$ decay exponentially and dynamically reduce to the classical transverse-traceless regime.

The results obtained in the present paper demonstrate that the Chen–Nagano gauge provides a natural geometric framework in which the scalar and trace-free components of Einstein metric deformations remain spectrally coupled. This coupling mechanism leads simultaneously to rigidity phenomena, spectral restrictions for the Lichnerowicz Laplacian, and dynamical stability results for the linearized Ricci flow (see [CK, Chapter 3, §2]).

In this sense, the Chen–Nagano approach extends the classical Berger–Ebin transverse-traceless theory by replacing the complete decoupling of the scalar and divergence-free components with a nontrivial spectral interaction governed by the geometry of the Einstein manifold and the curvature operator of the second kind.

We hope that the spectral interaction mechanism developed in this work will be useful for further investigations of Einstein metrics, geometric flows, and harmonic structures on Riemannian manifolds.

**Data Availability Statement**

No datasets were generated or analyzed during the current study.

**Funding Declaration**
The author declares that no financial funds, grants or other support were used in the preparation of this manuscript.

**References**

[B] Besse, A.L.: *Einstein Manifolds*. Ergeb. Math. Grenzgeb. (3), vol. 10. Springer, Berlin (1987)

[C] Chen, B.-Y.: Harmonic metrics, harmonic tensors and their applications. In: *Riemannian Geometry and Applications: Proceedings RIGA 2021*, pp. 17–44. Bucharest University Press, Bucharest (2021)

[CGT] Cao X., Gursky M.J., Tran H., Curvature of the second kind and a conjecture of Nishikawa, Comment. Math. Helv. 98 (2023), 195–216.

[CK] Chow, B., Knopf, D., *The Ricci Flow: An Introduction*. Mathematical Surveys and Monographs, vol. 110. American Mathematical Society, Providence (2004).

[CN] Chen, B.-Y., Nagano, T.: Harmonic metrics, harmonic tensors, and Gauss maps, *J. Math. Soc. Japan* **36**(2), 295–313 (1984)

[FP] Fine J., Premoselli B., Examples of compact Einstein four-manifolds with negative curvature, *Journal of the American Mathematical Society*, 33(4), (2020), 991–1038.

[J] Jäckel F., *Stability and construction of negatively curved Einstein metrics*, Ph.D. thesis, Rheinische Friedrich-Wilhelms-Universität Bonn, 2025.

[K] Kashiwada, T., *On the Curvature Operator of the Second Kind*. Natural Science Report, Ochanomizu University, 44(2) (1993), 69–73.

[K1] Koiso, N., Rigidity and stability of Einstein metrics. *Osaka J. Math.* **15**(2), 419–433 (1978)

[K2] Koiso, N.: Non-deformability of Einstein metrics. *Osaka J. Math.* **17**(2), 441–460 (1980)

[MRS] Mikeš J., Rovenski V., Stepanov S., An example of Lichnerowicz-type Laplacian, Annals of Global Analysis and Geometry, 58, (2020), 19–34.